\documentclass[11pt]{article}

\usepackage[
  left=1.2in,
  right=1.2in,
  top=0.8in,
  bottom=0.8in
]{geometry}

\usepackage{amsmath,amssymb,amsthm,mathtools}
\usepackage{algpseudocode}
\usepackage{bm}
\usepackage{booktabs}
\usepackage{graphicx}
\usepackage{array}
\usepackage{siunitx}
\usepackage{enumitem}
\usepackage{indentfirst}
\usepackage[numbers,sort&compress]{natbib}
\usepackage{xcolor}
\usepackage[colorlinks=true,linkcolor=blue,citecolor=blue,urlcolor=blue]{hyperref}
\usepackage{algorithm}

\newtheorem{remark}{Remark}

\graphicspath{
  {figures/}
}

\usepackage{caption}

\graphicspath{{figures/TV_ANM_result/}}

\begin{document}
\noindent\begin{minipage}{0.99\textwidth}
\centering
{\LARGE Adaptive Nesterov Momentum Method for Electrical Impedance Tomography with the Complete Electrode Model\par}
\vspace{1.1em}

{\large Kai Zhu$^{1}$, Jijun Liu$^{1}$, Min Zhong$^{1*}$\par}
\vspace{0.7em}

{\normalsize
$^{1}$Department of Mathematics, Southeast University,\\
Nanjing 210096, Jiangsu, P. R. China.\par}
\vspace{0.7em}

{\normalsize
*Corresponding author. E-mail:
\href{mailto:min.zhong@seu.edu.cn}{min.zhong@seu.edu.cn}.\par}
\end{minipage}

\begin{abstract}
We apply the adaptive Nesterov momentum (ANM) method~\cite{Jin2025AdaptiveNesterov} to electrical impedance tomography under the complete electrode model. 
The forward problem is formulated in the variational CEM setting, accounting for finite electrode size, contact impedance, insulating gaps, and the mean-free voltage gauge.
The resulting nonlinear inverse problem is treated within a unified dual-to-primal framework using three classes of strongly convex structural penalties: an \(L^2\)-type penalty, an \(L^1\)-type penalty promoting sparse deviations from a calibrated homogeneous background, and a TV-type penalty favoring approximately piecewise-constant conductivities with sharp interfaces. 
The TV class is implemented using both smoothed TV and Huber-TV formulations.
The data-misfit gradient is computed through CEM adjoint equations and stabilized by Sobolev smoothing. 
The method is evaluated on the publicly available KIT4 tank measurement data after calibration of the background conductivity and contact impedance from no-object measurements. 
The experiments include single, multiple, mixed-conductivity, and geometrically challenging phantom configurations. 
The \(L^2\)-type penalty generally produces smooth but diffuse reconstructions, whereas the \(L^1\)-type penalty yields cleaner backgrounds with occasional geometric distortion. 
The smoothed TV and Huber-TV penalties provide more spatially coherent localization and exhibit similar reconstruction behavior across most tested configurations. 
These results demonstrate the practical applicability of the adaptive Nesterov framework to measured CEM-EIT data.
\end{abstract}

\noindent\textbf{Keywords.}
electrical impedance tomography; complete electrode model; adaptive Nesterov momentum method; measured KIT4 data

\section{Introduction}
\label{sec:intro}
Electrical impedance tomography (EIT) is a non-invasive imaging technique that aims to recover the interior conductivity distribution of a physical body from boundary current--voltage measurements \cite{EIT99}. 
In a typical EIT experiment, electrical currents are applied through boundary electrodes, and the resulting electrode voltages are measured to reconstruct the internal conductivity distribution.
Owing to its low cost, portability, and absence of ionizing radiation, EIT has attracted considerable interest in medical imaging, industrial process monitoring, and nondestructive testing; see, for example, \cite{EIT_Andy, Holder2004EIT, Borcea2002EIT}.

From the mathematical viewpoint, conductivity reconstruction in EIT is a nonlinear and severely ill-posed inverse boundary value problem, commonly referred to as Calder\'on's inverse conductivity problem \cite{Calder2006EIT}. 
Its theoretical foundations have been studied extensively. Fundamental global uniqueness results were established under suitable regularity assumptions on the conductivity \cite{uniqueness1987EIT, Aless1990EIT}. 
Quantitative stability estimates and extensions to less regular and discontinuous conductivities were subsequently developed in \cite{Aless1988EIT,stability2001EIT,stability2010EIT}.
For detailed reviews of the mathematical theory of EIT, we refer the reader to \cite{Uhlmann2009EIT, EIT99}.

The classical Calder\'on formulation assumes idealized boundary data, whereas practical measurements are obtained through a finite number of electrodes.
For this setting, the complete electrode model (CEM) provides a realistic description by accounting for the finite electrode size, insulating gaps, electrode--medium contact impedances, and the grounding of the electrode voltages \cite{EIT_SIAM_92, Hyv2004CEM}. 
The influence of imperfectly known electrode positions, contact impedances, and other modeling parameters has also been investigated in \cite{Hyv2012electrode, Hyv2014electrode, winkler2016model}. 
The existence and uniqueness of solutions to the CEM forward problem were established in \cite{EIT_SIAM_92}, while its approximation properties and related inverse characterization results were further studied in \cite{Hyv2004CEM}. 
For conductivities that are piecewise constant on a prescribed pixel partition, uniqueness and Lipschitz stability from finitely many electrode measurements were established in \cite{Bastian2019EIT}.
Related results have shown that Calder\'on's problem with finitely many unknown parameters admits an equivalent convex semidefinite formulation, including extensions to finite measurement data \cite{Bastian2023SIAM_EIT}.

EIT reconstruction methods are commonly classified into difference and absolute imaging. 
Difference imaging seeks conductivity changes relative to a reference state and often benefits from the partial cancellation of systematic measurement and modeling errors. 
Computationally inexpensive reconstructions can be obtained by linearizing the forward map around a reference conductivity, as in NOSER and related backprojection or regularized linearized schemes \cite{Cheney1990NOSER, backprojection1990EIT, Amm2016linearEIT}.
Absolute imaging, by contrast, seeks to recover the conductivity itself from the measured electrode voltages and therefore requires the solution of a nonlinear inverse problem. 
Methods developed for this setting include geometric and level-set approaches \cite{levelset2005EIT, Alghamdi2024Spatial}, convex formulations for finite-dimensional conductivity classes \cite{Bastian2023SIAM_EIT}, Newton-type and inexact Newton regularization methods \cite{Rieder2006EIT, EIT_local_con_2008,winkler2016model}, and Levenberg--Marquardt-type methods \cite{rrLM_banach_IP_2022}.
Alternative reconstruction paradigms include D-bar and Calder\'on-type methods \cite{Dbar2000EIT, Muller2008EIT, Hamilton2018Dbar, Santos2020DbarPrior}, as well as direct regularized reconstruction methods for the three-dimensional Calder\'on problem \cite{direct2022EIT}. 
More recently, data-driven and deep-learning-based reconstruction methods have also been investigated \cite{Denker2024DataDrivenEIT, Denker2025DL}.

Structural prior information can be incorporated through quadratic smoothness, total variation (TV), and sparsity-promoting penalties. 
TV is particularly suitable for approximately piecewise-constant conductivities with sharp interfaces \cite{Gonzalez2016TV3D}, whereas \(L^1\)-type penalties promote sparse or localized deviations from a reference conductivity \cite{EIT_sparsity_2012}.

Despite these developments, stable and computationally efficient nonlinear reconstruction from noisy electrode data remains challenging. 
Landweber-type iterative regularization offers a flexible framework because it requires only evaluations of the forward operator and its adjoint derivative and permits nonsmooth structural information to be incorporated through a uniformly convex penalty functional \cite{kaltenbacher2008iterative, Jin2025}. 
Its main limitation is the typically slow convergence of the basic iteration. 
Motivated by recent adaptive Nesterov momentum strategies for ill-posed inverse problems \cite{Humber_Ramlau_2017_TPG, Jin2025AdaptiveNesterov, zz2025}, we employ an adaptive Nesterov momentum method for the CEM-EIT inverse problem. 
The method updates a dual variable through an adjoint-gradient step and a momentum extrapolation, while the conductivity is recovered by a dual-to-primal minimization.

The main contribution of this work is the application and measured-data evaluation of the adaptive Nesterov momentum method for the nonlinear CEM-EIT inverse problem with \(L^2\)-type, \(L^1\)-type, and TV-type structural penalties, where the TV class is implemented using both smoothed TV and Huber-TV formulations. 
The numerical implementation is adapted to the computational structure of the complete electrode model. 
To improve numerical robustness, Sobolev smoothing is applied to the raw adjoint gradient to suppress oscillatory update directions and stabilize the reconstruction.

The resulting reconstruction procedure is evaluated on measured data from the public KIT4 tank data set \cite{hauptmann2018eitdataset}. 
The experiments show that the adaptive Nesterov momentum method can be used effectively for measured CEM-EIT data in the presence of measurement noise and modeling discrepancies. 
Within this setting, the \(L^2\)-type penalty provides a smooth baseline but tends to produce more spatially diffuse conductivity variations. 
The \(L^1\)-type penalty promotes sparse deviations from a calibrated homogeneous background and admits an inexpensive dual-to-primal update, whereas the TV-type penalty favors piecewise-constant conductivity distributions with sharper and more spatially coherent interfaces.

The remainder of the paper is organized as follows.
Section~\ref{sec:cem} introduces the complete electrode model (CEM), together with its variational formulation, the finite-data forward operator, and the associated derivative and adjoint problems.
Section~\ref{sec:anm} develops the adaptive Nesterov momentum framework, including the dual-to-primal formulation, the incorporation of structural penalties, and the Sobolev-smoothed adjoint gradient.
Section~\ref{sec:numerics} presents numerical reconstructions from measured KIT4 data and compares the performance of \(L^2\)-type, \(L^1\)-type, and TV-type penalties, with particular attention to the smoothed-TV and Huber-TV formulations.
Additional details on the selection of the strong-convexity parameter are provided in Appendix~\ref{app:beta-selection}.

\section{EIT with the Complete Electrode Model}
\label{sec:cem}

We recall the complete electrode model (CEM) used throughout the paper. 
Let \(\Omega\subset\mathbb R^d\), \(d=2,3\), be a bounded conductive domain. 
The conductivity is assumed to satisfy 
$
    \sigma\in L^\infty_+(\Omega)
    :=
    \left\{
    \kappa\in L^\infty(\Omega):
    \operatorname*{ess\,inf}_{x\in\Omega}\kappa(x)>0
    \right\}.
$
Let \(e_\ell\subset\partial\Omega\), \(\ell=1,\ldots,L\), be open, mutually disjoint electrodes with positive surface measure, and let \(z_\ell>0\) denote the contact impedance on \(e_\ell\). 
For a charge-conserving current pattern
\[
    I=(I_1,\ldots,I_L)^T\in\mathbb R^L,
    \qquad
    \sum_{\ell=1}^L I_\ell=0,
\]
the CEM determines the interior electric potential \(u\) and the electrode voltages \(U=(U_1,\ldots,U_L)^T\) from
\begin{equation}
\label{eq:cem-strong}
\left\{
\begin{aligned}
-\nabla\cdot(\sigma\nabla u) &= 0
&& \text{in } \Omega,\\
u+z_\ell\sigma\frac{\partial u}{\partial \nu} &= U_\ell
&& \text{on } e_\ell,\quad \ell=1,\ldots,L,\\
\int_{e_\ell}\sigma\frac{\partial u}{\partial \nu}\,ds &= I_\ell
&& \ell=1,\ldots,L,\\
\sigma\frac{\partial u}{\partial \nu} &= 0
&& \text{on } \partial\Omega\setminus\bigcup_{\ell=1}^L e_\ell .
\end{aligned}
\right.
\end{equation}
Here \(\nu\) denotes the unit outward normal vector. 
Since the electric potential is determined only up to an additive constant, we impose the grounding condition $\sum_{\ell=1}^L U_\ell=0$.
We use the notation
\[
    \mathbb R^L_\diamond
    :=
    \left\{
    V\in\mathbb R^L:
    \sum_{\ell=1}^L V_\ell=0
    \right\},
    \qquad
    \mathcal H:=H^1(\Omega)\oplus\mathbb R^L_\diamond .
\]
The corresponding variational formulation is: find \((u,U)\in\mathcal H\) such that
\begin{equation}
\label{eq:var-CEM-EIT}
    \mathcal B_{\sigma,z}((u,U),(v,V))
    =
    \sum_{\ell=1}^L I_\ell V_\ell,
    \qquad
    \forall (v, V)\in\mathcal H,
\end{equation}
where
\[
\mathcal B_{\sigma,z}((u,U),(v,V))
:=
\int_\Omega \sigma\nabla u\cdot\nabla v\,dx
+
\sum_{\ell=1}^L
\frac{1}{z_\ell}
\int_{e_\ell}
(u-U_\ell)(v-V_\ell)\,ds .
\]
For each \(I\in\mathbb R^L_\diamond\), the CEM is well posed in \(\mathcal H\); see \cite{EIT_SIAM_92}.

We now introduce the finite-measurement forward map. 
Let \(I^{(i)}\in\mathbb R^L_\diamond\), \(i=1,\ldots,N\), be the prescribed current patterns, and let \((u^{(i)},U^{(i)})\in\mathcal H\) solve \eqref{eq:var-CEM-EIT}. 
The CEM current-to-voltage map is
\[
    \Lambda_{\sigma,z}:\mathbb R^L_\diamond\to\mathbb R^L_\diamond,
    \qquad
    \Lambda_{\sigma,z}I^{(i)}=U^{(i)} .
\]
The contact impedance $z$ is estimated beforehand and fixed during the conductivity reconstruction.
We therefore write
\begin{equation}
\label{eitforward}
    F_i(\sigma):=\Lambda_{\sigma,z}I^{(i)} =U^{(i)},
    \qquad i=1,\ldots,N,
\end{equation}
and define the stacked forward map
\begin{equation}
\label{eq:stacked-forward-map}
    F(\sigma)
    :=
    \big(F_1(\sigma),\ldots,F_N(\sigma)\big)
    \in Y,
    \qquad
    Y:=(\mathbb R^L_\diamond)^N .
\end{equation}
Given noisy electrode-voltage measurements
\[
    U^\delta
    =
    \bigl(U^{(1),\delta},\ldots,U^{(N),\delta}\bigr)\in Y,
\]
where the superscript \(\delta\) indicates that the measured data may be
contaminated by noise, we assume that $\|U-U^\delta\| \leq \delta$,
where \(U\) denotes the exact electrode-voltage data.
The inverse problem is then to recover \(\sigma\) from
\[
    F(\sigma)\approx U^\delta .
\]
This finite-measurement inverse problem is nonlinear and ill-posed and must therefore be regularized. 
We employ an iterative regularization method in which the conductivity is updated through residual-based gradient corrections in the dual space. 
The evaluation of these corrections requires the Fr\'echet derivative of the CEM forward map and the action of its adjoint.
We next recall the derivative formula and derive the corresponding adjoint
representation.

For
\(\sigma\in\operatorname{int}(L^\infty_+(\Omega))\) and
\(h\in L^\infty(\Omega)\), the derivative of \(F_i\) in the direction \(h\) is given by
\[
    F_i'(\sigma)h=W^{(i)},
\]
where \((w^{(i)},W^{(i)})\in\mathcal H\) solves
\begin{equation}
\label{eq:CEM-derivative}
    \mathcal B_{\sigma,z}
    \big((w^{(i)},W^{(i)}),(v,V)\big)
    =
    -\int_\Omega
    h\,\nabla u^{(i)}\cdot\nabla v\,dx,
    \qquad
    \forall (v, V)\in\mathcal H .
\end{equation}
Here \((u^{(i)},U^{(i)})\) is the CEM solution corresponding to \(I^{(i)}\). 
This derivative formula is standard; see \cite{Rieder2006EIT, winkler2016model}.
For \(Z\in\mathbb R^L_\diamond\), let \((p^{(i)},P^{(i)})\in\mathcal H\) solve the adjoint CEM problem
\begin{equation}
\label{eq:CEM-adjoint-state}
    \mathcal B_{\sigma,z}
    \big((p^{(i)},P^{(i)}),(v,V)\big)
    =
    \sum_{\ell=1}^L Z_\ell V_\ell,
    \qquad
    \forall (v, V)\in\mathcal H .
\end{equation}
Then the \(L^2(\Omega)\)-gradient density associated with the adjoint action is
\begin{equation}
\label{eq:CEM-adjoint-density}
    F_i'(\sigma)^*Z
    =
    -\nabla u^{(i)}\cdot\nabla p^{(i)} .
\end{equation}


\section{Adaptive Nesterov Momentum Method}
\label{sec:anm}

We reconstruct the conductivity $\sigma$ from noisy electrode-voltage data $U^\delta$ using the adaptive Nesterov momentum method proposed in \cite{Jin2025AdaptiveNesterov}. 
The method combines iterative regularization with a momentum extrapolation strategy. 
In contrast to a standard Landweber-type iteration, the gradient correction is performed from an extrapolated dual variable, so that information from two consecutive iterations is used to accelerate the reconstruction. 
The step size and the momentum parameter are selected adaptively through explicit formulas rather than by a line search. 
Their construction is coupled with the ordering of the gradient and extrapolation steps and is designed to retain the descent properties required by the convergence analysis.

Another important feature is that the conductivity is recovered from the dual variable through a strongly convex minimization problem. 
This dual-to-primal formulation allows different structural penalties and pointwise conductivity constraints to be incorporated. 
Consequently, the \(L^2\)-, \(L^1\)-, and TV-type reconstructions considered below differ only in the evaluation of the dual-to-primal map.

\subsection{Adaptive Nesterov Momentum Method}

Let $F$ denote the CEM parameter-to-data map introduced in Section~\ref{sec:cem}. 
Starting from the initial dual variables $\xi_0^\delta=\zeta_0^\delta$, the iteration is defined by
\begin{equation}
\label{eq:anm-update}
\begin{aligned}
\xi_{n+1}^\delta
&=
\zeta_n^\delta
-
\mu_n^\delta
F'(\sigma_n^\delta)^*
\bigl(F(\sigma_n^\delta)-U^\delta\bigr),
\\[0.3em]
\zeta_{n+1}^\delta
&=
\xi_{n+1}^\delta
+
\lambda_n^\delta
\bigl(\xi_{n+1}^\delta-\xi_n^\delta\bigr),
\\[0.3em]
\sigma_{n+1}^\delta
&=
\operatorname*{arg\,min}_{\sigma\in L^2(\Omega)}
\left\{
\Theta(\sigma)
-
\left\langle
\zeta_{n+1}^\delta,\sigma
\right\rangle
\right\}.
\end{aligned}
\end{equation}
Here, \(F'(\sigma_n^\delta)^*\) denotes the adjoint of the Fr\'echet derivative of the CEM forward map, \(\mu_n^\delta\) is the adaptive step size, and \(\lambda_n^\delta\) is the adaptive momentum parameter. 
Thus, one iteration consists of an adjoint-gradient correction, a Nesterov momentum extrapolation, and a dual-to-primal reconstruction.
The residual at iteration \(n\) is denoted by
$r_n^\delta
=
F(\sigma_n^\delta)-U^\delta$.

\begin{remark}
We conclude the description of the method with two observations concerning
its limiting case and computational cost.
\begin{itemize}
    \item
    When the momentum parameter vanishes, that is,
    \(\lambda_n^\delta=0\), the scheme reduces to the corresponding
    Landweber-type iteration.

    \item
    The main computational cost of each iteration arises from the forward
    and adjoint CEM solves. The explicit updates of \(\mu_n^\delta\) and
    \(\lambda_n^\delta\), together with the \(L^2\)- and \(L^1\)-type
    dual-to-primal maps, have negligible computational cost. By contrast,
    the TV-type dual-to-primal update requires the numerical solution of a
    constrained denoising subproblem and may therefore contribute
    substantially to the total computational cost; see
    Subsection~\ref{subsec:anm-penalties}.
\end{itemize}
\end{remark}

\begin{algorithm}[H]
\caption{Adaptive Nesterov momentum method}
\label{alg:adaptive-nesterov}
\begin{algorithmic}[1]

\State Choose an initial dual variable $\zeta_0\in L^2(\Omega)$ and compute
\Statex $
\sigma_0 := \operatorname*{arg\,min}_{\sigma\in L^2(\Omega)}
\left\{
\Theta(\sigma) - \langle \zeta_0,\sigma\rangle
\right\}.
$

\State Set $
\sigma_{-1}^{\delta} = \sigma_0^{\delta} = \sigma_0,
\quad \xi_0^{\delta} = \zeta_0^{\delta} = \zeta_0,\quad \lambda_{-1}^{\delta}=0,
\quad \widetilde{\gamma}_0^{\delta}=0, \quad m_0^{\delta}=0.
$

\State Choose $\tau>1$ such that
$c_{\eta,\tau}=1-\eta-(1+\eta)/\tau>0$,
together with $0<\mu_0$, $\mu_1>0$, a nonnegative sequence
$\{\hat{\lambda}_k\}_{k\ge0}$, and the strong convexity constant
$\kappa_\Theta>0$ of $\Theta$.

\For{$k=0,1,2,\ldots$}

    \State Compute the residual $ r_k^{\delta}   := F(\sigma_k^{\delta})-U^{\delta}. $

    \If{$\|r_k^{\delta}\|\leq\tau\delta$}
        \State Set $n_{\delta}:=k$.
        \State \textbf{break}
    \EndIf

    \State Compute the dual-gradient direction
$
    g_k^{\delta} := F'(\sigma_k^{\delta})^*(r_k^{\delta}).
$

    \State Choose the step size
$
    \mu_k^{\delta}
    :=
    \min
    \left\{
    \mu_0
    \frac{
    c_{\eta,\tau}
    \|r_k^{\delta}\|^2
    }{
    \|g_k^{\delta}\|^2
    },
    \;
    \mu_1
    \right\}.
$

    \State Compute the auxiliary quantity
    \[
    \begin{aligned}
    \widetilde{\gamma}_{k+1}^{\delta}
    :=
    \lambda_{k-1}^{\delta}
    \left\langle
    m_k^{\delta},
    \sigma_k^{\delta}-\sigma_{k-1}^{\delta}
    \right\rangle
    +
    \lambda_{k-1}^{\delta}
    \widetilde{\gamma}_k^{\delta}- c_{\eta,\tau}\mu_k^{\delta}
    \|r_k^{\delta}\|^2.
    \end{aligned}
    \]

    \State Update the intermediate dual variable
    $
    \xi_{k+1}^{\delta}:=\zeta_k^{\delta}-\mu_k^{\delta}g_k^{\delta}.
    $

    \State Compute the momentum direction
    $
    m_{k+1}^{\delta}
    :=
    \xi_{k+1}^{\delta}
    -
    \xi_k^{\delta}.
    $
    \State Compute the adaptive momentum parameter $\lambda_k^\delta$ by
    \[
    \lambda_k^\delta=
    \begin{cases}
    0,
    & m_{k+1}^\delta=0,\\[0.8ex]
    \displaystyle
    \min\left\{
    \max\left\{
    0,\,
    \frac{\mu_k^\delta\left\langle g_k^\delta,m_{k+1}^\delta\right\rangle - 2\kappa_\Theta\,\widetilde{\gamma}_{k+1}^\delta}
    {\|m_{k+1}^\delta\|^2}
    \right\},
    \hat{\lambda}_k
    \right\},
    &
    m_{k+1}^\delta\neq0.
    \end{cases}
    \]

    \State Update the accelerated dual variable
    $
    \zeta_{k+1}^{\delta}:=\xi_{k+1}^{\delta} +\lambda_k^{\delta}m_{k+1}^{\delta}.
    $

    \State Update the primal iterate
    $
    \sigma_{k+1}^{\delta}
    :=
    \operatorname*{arg\,min}_{\sigma\in L^2(\Omega)}
    \left\{
    \Theta(\sigma)
    -
    \left\langle
    \zeta_{k+1}^{\delta},
    \sigma
    \right\rangle
    \right\}.
    $

\EndFor

\State \textbf{Output:}
the discrepancy-stopped reconstruction
$\sigma_{n_{\delta}}^{\delta}$.

\end{algorithmic}
\end{algorithm}

\begin{remark}
The algorithmic parameters are defined as follows.
\begin{itemize}
    \item The parameter $\tau$ is chosen such that
    $c_{\eta,\tau}>0$, where $\eta$ is the tangential cone constant associated
    with the nonlinear forward operator.

    \item The sequence $\{\hat{\lambda}_k\}_{k\ge0}$ specifies the
    Nesterov extrapolation. A typical choice is
    \[
        \hat{\lambda}_k=\frac{k}{k+\alpha},
        \qquad \alpha\geq3,
    \]
    which ensures $0\leq\lambda_k^\delta\leq\hat{\lambda}_k<1$ at every finite iteration.
    \item Under the strong convexity convention in
    \cite{Jin2025AdaptiveNesterov}, $\kappa_\Theta=1/2$ for
    $\Theta_{L^2}$ and $\kappa_\Theta=1/(2\beta)$ for
    $\Theta_{L^1}$ and $\Theta_{\rm TV}$.
\end{itemize}
\end{remark}

In the numerical implementation, the adjoint-gradient direction is smoothed
in a Sobolev metric before being used in the dual update. For
\(r=(r^{(1)},\ldots,r^{(N)})\in Y\), the stacked raw adjoint gradient is
\begin{equation}
\label{eq:CEM-stacked-gradient}
F'(\sigma)^*r
=
\sum_{i=1}^N F_i'(\sigma)^*r^{(i)}
=
-\sum_{i=1}^N
\nabla u^{(i)}\cdot\nabla p^{(i)},
\end{equation}
where \(p^{(i)}\) solves \eqref{eq:CEM-adjoint-state} with
\(Z=r^{(i)}\). We replace
\(g_{\rm raw}=F'(\sigma)^*r\) by its Sobolev-smoothed counterpart
\(g=\mathcal S_q g_{\rm raw}\), defined by
\begin{equation}
\label{sobolev-smooth}
\left\{
\begin{aligned}
-q\Delta g+g &= g_{\rm raw}
&&\text{in }\Omega,\\
g &=0
&&\text{on }\partial\Omega,
\end{aligned}
\right.,
\end{equation}
where \(q>0\) is the smoothing parameter. Accordingly, the
adjoint-gradient term in \eqref{eq:anm-update} is understood as its
Sobolev-smoothed version. This step is intended to suppress small-scale
oscillations in the update direction; related Sobolev-gradient techniques
have been used in EIT reconstruction~\cite{EIT_sparsity_2012}.

\subsection{Structural penalty functionals}
\label{subsec:anm-penalties}
The functional \(\Theta\) is chosen to satisfy the strong convexity requirement of the adaptive Nesterov momentum method while incorporating structural prior information about the conductivity. 
The indicator term further enforces the admissible conductivity constraints.

Let
\[
\mathcal A
=
\left\{
\sigma\in L^\infty(\Omega):
c\leq\sigma\leq c^{-1}\ \text{a.e. in }\Omega
\right\},
\]
where \(c\in(0,1)\) is fixed, and let
\[
\iota_{\mathcal A}(\sigma)
=
\begin{cases}
0, & \sigma\in\mathcal A,\\
+\infty, & \sigma\notin\mathcal A,
\end{cases}
\]
denote the indicator functional of the admissible set. The bounds encoded in
\(\mathcal A\) enforce positivity of the conductivity and prevent
nonphysical iterates.

We consider three classes of functionals \(\Theta\). The \(L^2\)-type
functional provides a smooth quadratic baseline, the \(L^1\)-type functional
promotes sparse deviations from a homogeneous background, and the TV-type
functional favors approximately piecewise-constant conductivity distributions
with sharp interfaces.
\paragraph{\(L^2\)-type penalty}

As a smooth baseline, we consider
\begin{equation}
\label{eq:l2-penalty}
    \Theta_{L^2}(\sigma)
    =
    \frac{1}{2}
    \|\sigma\|_{L^2(\Omega)}^2
    +
    \iota_{\mathcal A}(\sigma).
\end{equation}
The corresponding dual-to-primal update is explicit and reduces to a
pointwise projection onto the admissible interval.
\paragraph{$L^1$-type penalty}

Let $\sigma_{\rm bk}$ denote the calibrated homogeneous background conductivity. 
The background-centered $L^1$-type functional is defined by
\begin{equation}
\label{eq:l1-penalty}
\Theta_{L^1}(\sigma)
=
\frac{1}{2\beta}\|\sigma\|_{L^2(\Omega)}^2
+
\|\sigma-\sigma_{\rm bk}\|_{L^1(\Omega)}
+
\iota_{\mathcal A}(\sigma),
\end{equation}
where $\beta>0$. 
The quadratic term makes the functional strongly convex, whereas the $L^1$ term promotes spatially localized deviations from the homogeneous background.

Because the nonsmooth term in \eqref{eq:l1-penalty} is pointwise separable, the dual-to-primal problem admits an explicit solution. Substituting $\Theta_{L^1}$ into the last equation of \eqref{eq:anm-update} gives
\begin{equation}
\label{eq:l1-prox}
\sigma_{n+1}^\delta
=
P_{\mathcal A}
\left[
\sigma_{\rm bk}
+
S_\beta\!\left(
\beta\zeta_{n+1}^\delta-\sigma_{\rm bk}
\right)
\right],
\end{equation}
where
\[
S_\beta(t)=\operatorname{sign}(t)\max\{|t|-\beta,0\}
\]
is the soft-thresholding operator, and \(P_{\mathcal A}\) denotes the pointwise projection onto the admissible interval.

\paragraph{Total variation}

For conductivity distributions that are approximately piecewise constant, we consider the total variation functional
\begin{equation}
\label{eq:tv-penalty}
\Theta_{\rm TV}(\sigma)
=
\frac{1}{2\beta}\|\sigma\|_{L^2(\Omega)}^2
+
{\rm TV}(\sigma)
+
\iota_{\mathcal A}(\sigma),
\end{equation}
where
\begin{equation}
\label{eq:tv-definition}
{\rm TV}(\sigma)
=
\int_\Omega |\nabla \sigma|\,{\rm d}x,
\end{equation}
for sufficiently regular \(\sigma\), with the standard extension to \(BV(\Omega)\). 
The TV seminorm suppresses small-scale variations while allowing sharp transitions across material interfaces.

Unlike the pointwise \(L^1\) penalty, the TV functional couples neighboring spatial degrees of freedom through the conductivity gradient. 
Consequently, the associated dual-to-primal map is not pointwise separable and is given by
\begin{equation}
\label{eq:tv-prox}
\sigma_{n+1}^\delta
=
\operatorname*{arg\,min}_{\sigma\in\mathcal A}
\left\{
\frac{1}{2\beta}
\left\|
\sigma-\beta\zeta_{n+1}^\delta
\right\|_{L^2(\Omega)}^2
+
{\rm TV}(\sigma)
\right\}.
\end{equation}
Thus, each TV-based dual-to-primal update requires the numerical solution of a box-constrained Rudin--Osher--Fatemi (ROF)-type denoising problem.

\subparagraph{Smoothed TV}

Since the standard TV density is nondifferentiable at \(\nabla\sigma=0\), we employ the differentiable approximation
\begin{equation}
\label{eq:smoothed-tv}
{\rm TV}_{\varepsilon}(\sigma)
=
\int_\Omega
\sqrt{|\nabla\sigma|^2+\varepsilon}
\,{\rm d}x,
\qquad
\varepsilon>0.
\end{equation}
The parameter \(\varepsilon\) removes the singularity at \(\nabla\sigma=0\) and permits a Newton-type treatment. 
The resulting dual-to-primal subproblem is solved inexactly by a primal--dual Newton method, abbreviated as PD--NT, applied to the corresponding coupled primal--dual optimality system~\cite{PDnewtonmethod}.

\subparagraph{Huber-TV}

As an alternative regularized TV formulation, we consider
the Huber-TV functional~\cite{AscherHaberHuang2006}
\begin{equation}
\label{eq:huber-tv}
{\rm TV}_{\gamma}(\sigma)
=
\int_\Omega
\psi_\gamma(|\nabla\sigma|)
\,{\rm d}x,
\end{equation}
where the normalized Huber density~\cite{Huber1964} is defined by
\begin{equation}
\label{eq:huber-function}
\psi_\gamma(t)
=
\begin{cases}
\dfrac{t^2}{2\gamma},
& 0\leq t\leq\gamma,\\[1ex]
t-\dfrac{\gamma}{2},
& t>\gamma,
\end{cases}
\qquad
\gamma>0.
\end{equation}
The Huber density is quadratic below the threshold $\gamma$
and agrees with the standard TV density, up to an additive
constant, above this threshold. It is continuously differentiable
but not twice differentiable at $t=\gamma$. Thus, it smooths
small gradients while retaining the linear growth associated
with the edge-preserving behavior of TV for larger gradients.

Both TV formulations lead to spatially coupled dual-to-primal
subproblems and are therefore more computationally involved
than the explicit $L^1$-type mapping. Their reconstruction
behavior is examined in Section~\ref{sec:numerics}.

\section{Validation on Measured KIT4 Data}
\label{sec:numerics}

We evaluate the adaptive Nesterov momentum method using the publicly available KIT4 tank measurement data set \cite{hauptmann2018eitdataset}.
Three classes of structural penalties are considered within the same CEM-based adaptive Nesterov framework: an \(L^2\)-type penalty promoting smooth reconstructions, an \(L^1\)-type penalty promoting sparse deviations from a homogeneous background, and a TV-type penalty favoring piecewise-constant conductivities with sharp interfaces. 
The TV-type penalty is implemented using both smoothed TV and Huber-TV formulations. 
Their reconstruction quality and residual behavior are evaluated on several measured KIT4 tank configurations.

\subsection{KIT4 data set and computational details}

\begin{figure}[H]
    \centering

    \begin{minipage}[t]{0.46\textwidth}
        \centering
        \includegraphics[
            width=\linewidth,
            height=0.32\textheight,
            keepaspectratio
        ]{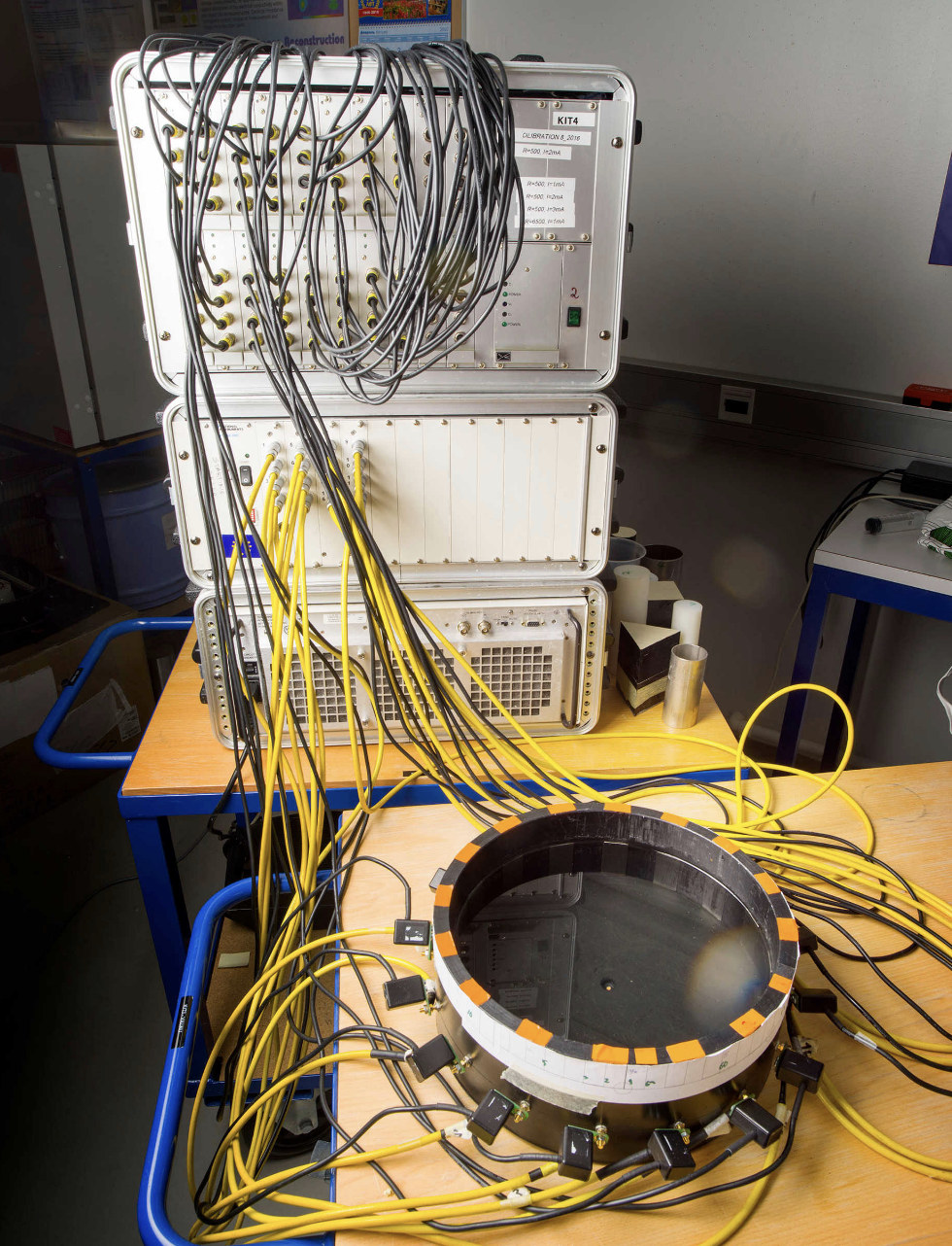}

        \smallskip
        \textbf{(a)} KIT4 EIT measurement system.
    \end{minipage}
    \hfill
    \begin{minipage}[t]{0.50\textwidth}
        \centering
        \includegraphics[
            width=\linewidth,
            height=0.32\textheight,
            keepaspectratio
        ]{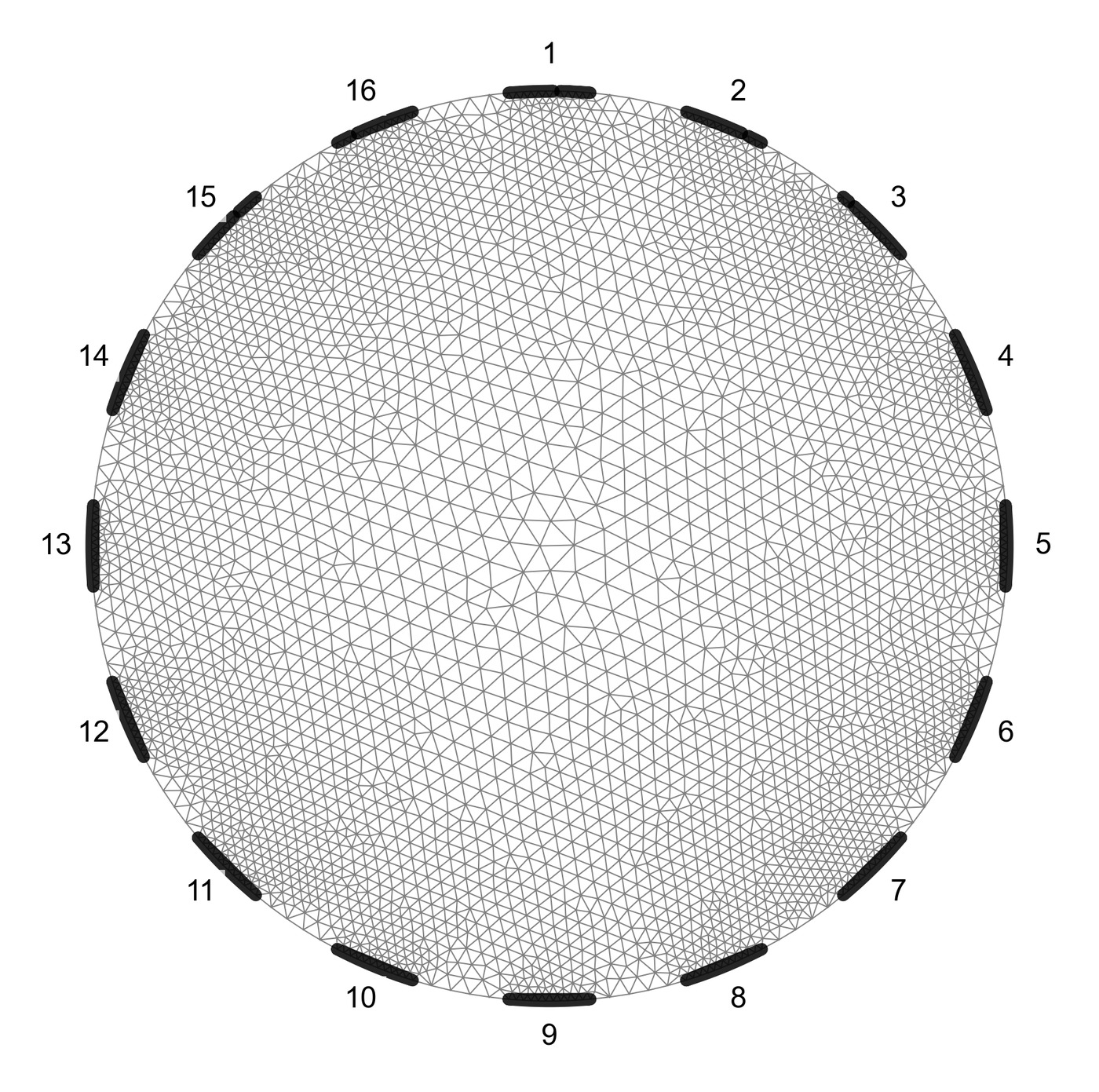}

        \smallskip
        \textbf{(b)} Finite element mesh and electrode configuration.
    \end{minipage}

    \caption{Experimental measurement system and computational model.
    The KIT4 data were acquired at the University of Eastern Finland;
    the system photograph is reproduced from
    \cite{hauptmann2018eitdataset}.}
    \label{fig:kit4-setup}
\end{figure}

The experiments are based on the publicly available KIT4 EIT tank data set \citep{hauptmann2018eitdataset}, acquired using the Kuopio Impedance Tomography (KIT4) system. 
The measurements were performed in a cylindrical saline tank of diameter $28\,\mathrm{cm}$ equipped with $16$ rectangular boundary electrodes, each of height $7\,\mathrm{cm}$ and width $2.5\,\mathrm{cm}$; see Fig.~\ref{fig:kit4-setup}(a).
The phantom configurations were generated by placing conductive or insulating targets of various shapes and arrangements in a tank filled with homogeneous saline.
We use $N=16$ adjacent current-injection patterns, each applying a current of amplitude $2\,\mathrm{mA}$ between neighboring electrodes.

Owing to the essentially two-dimensional geometry of the cylindrical tank and the vertically extended targets, the complete electrode model \eqref{eq:cem-strong} and the corresponding adjoint problem \eqref{eq:CEM-adjoint-state} were discretized using the finite element method.
The computational domain \(\Omega\) was triangulated into \(7,132\) triangular elements with \(3,695\) vertices; see Fig.~\ref{fig:kit4-setup}(b). 
The mesh was locally refined near the boundary electrodes to better resolve the electrode geometry and the associated electrode--medium interface terms.

Before reconstructing the inclusions, the background conductivity and contact impedance were calibrated from the no-object measurements. 
The calibrated background conductivity was then used in the \(L^1\)-type penalty, while the estimated contact impedance was kept fixed in the subsequent forward solves.
Assuming a homogeneous background conductivity \(\sigma_{\rm bk}\) and identical contact impedances \(z_\ell=z\), we used the CEM calibration procedure in \citep{winkler2016model}.
The measured voltage differences were first converted into mean-free electrode potentials, after which \(\sigma_{\rm bk}\) and \(z\) were estimated from the background data. 
The calibrated values used in all subsequent reconstructions were
\[
\sigma_{\rm bk}
=
1.8723\times10^{-3}\,\mathrm{S},
\qquad
z
=
2.5\times10^{-4} \Omega \cdot \text{cm}.
\]

\begin{figure}[H]
    \centering
    \includegraphics[
        width=0.75\textwidth,
        height=0.45\textheight,
        keepaspectratio
    ]{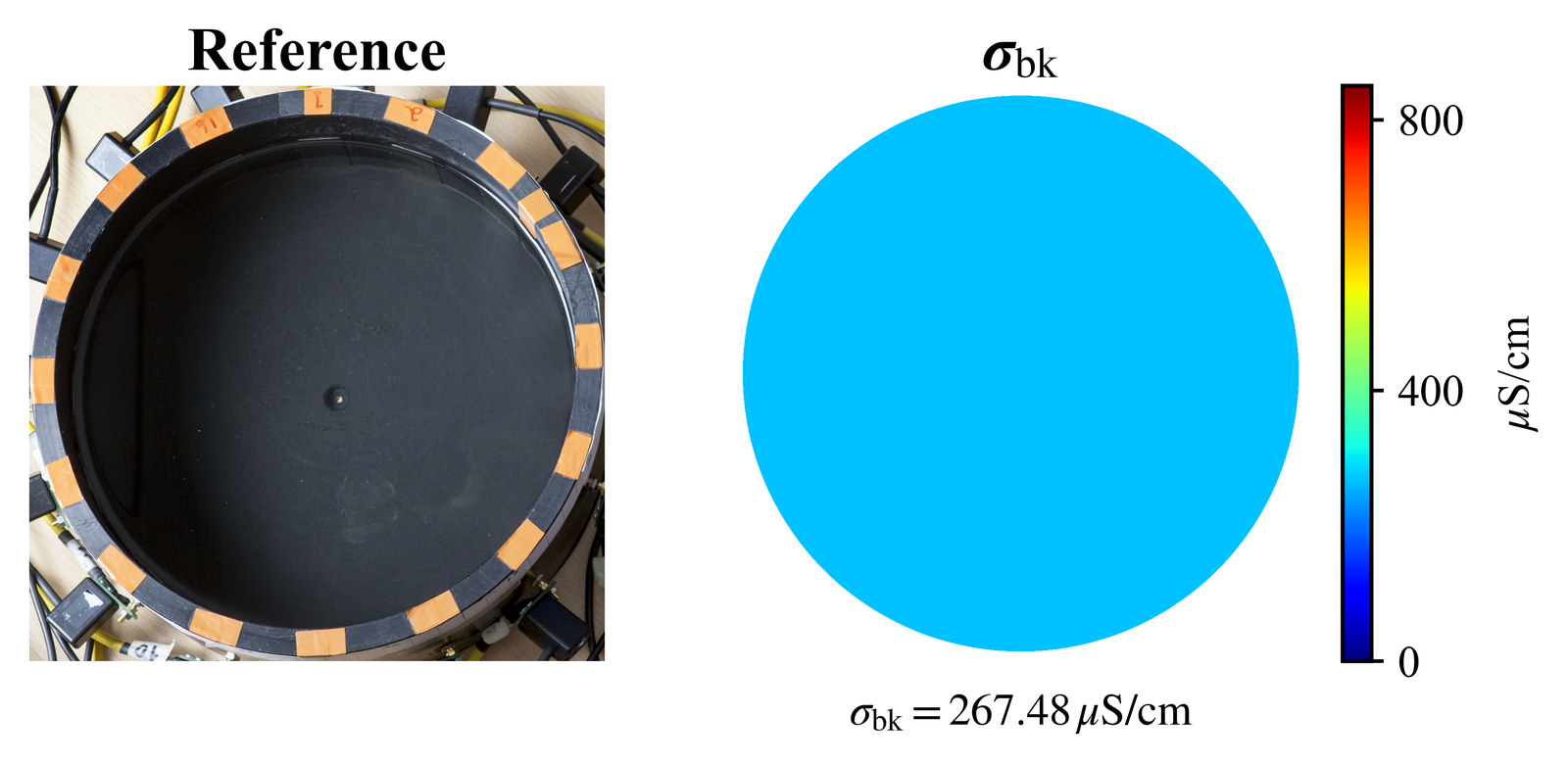}
    \caption{Background calibration using the no-object KIT4 measurements.
    The left panel shows the homogeneous saline tank, and the right panel shows
    the conductivity reconstruction obtained under the homogeneous
    single-parameter model.}
    \label{fig:background-calibration}
\end{figure}

The calibrated \(\sigma_{\rm bk}\) is the effective two-dimensional sheet conductivity. 
Using the effective tank height \(h=7\,\mathrm{cm}\), the corresponding three-dimensional bulk conductivity is
\[
\frac{\sigma_{\rm bk}}{h}
=
267.48\,\mu\mathrm{S/cm}.
\]
As shown in Fig.~\ref{fig:background-calibration}, the no-object tank contains only the homogeneous saline background, and the corresponding single-parameter reconstruction is nearly spatially uniform.
This value is reasonably consistent with the experimentally reported saline conductivity of \(300\,\mu\mathrm{S/cm}\) at \(19\,^{\circ}\mathrm{C}\).
All subsequent reconstruction results are converted to the corresponding three-dimensional bulk conductivity and displayed in \(\mu\mathrm{S/cm}\).

\subsection*{Numerical setting}
We briefly summarize the numerical parameters used in Algorithm~\ref{alg:adaptive-nesterov}. 
The initial dual variable, the step-size parameters, the prescribed upper bound for the momentum parameter, and the Sobolev smoothing parameter are chosen as
\[
    \zeta_0(x)\equiv 1,\qquad
    \mu_0=1.9\kappa_\Theta,\qquad
    \mu_1=600,\qquad
    \widehat{\lambda}_n=\frac{n}{n+3},\qquad
    q=10^{-2}.
\]
The tangential cone constant is estimated as \(\eta=0.25\), following \cite{winkler2016model}, while the noise level is estimated as \(\delta=0.0314\). 
The discrepancy parameter is set to \(\tau=1.75\), which satisfies the required condition \(c_{\eta,\tau}>0\). 
The maximum number of outer iterations is set to \(k_{\max}=800\).
Since no theoretical parameter-choice rule for \(\beta\) is available in the
present measured-data setting, the parameters are selected empirically through
the preliminary study reported in Appendix~\ref{app:beta-selection}. Based on
the balance among background suppression, target contrast, and geometric
fidelity, we use
\[
    \beta_{L^1}=5,\qquad
    \beta_{\mathrm{TV}_\varepsilon}=2,\qquad
    \beta_{\mathrm{TV}_\gamma}=2
\]
in all subsequent experiments. The smoothing parameters are fixed at
\(\varepsilon=10^{-6}\) and \(\gamma=10^{-3}\)
for the smoothed TV penalty \(\mathrm{TV}_\varepsilon\) and the Huber-TV penalty
\(\mathrm{TV}_\gamma\), respectively.

\subsection{Reconstruction results}
\label{subsec:reconstruction-results}

The reconstruction results are presented in Figs.~\ref{fig:single-inclusions}, \ref{fig:two-inclusions}, \ref{fig:three-conductive-inclusions}, \ref{fig:mixed-inclusions}, and \ref{fig:challenging-phantoms}.
In each figure, the five columns show, from left to right, the reference target configuration and the \(L^2\)-, \(L^1\)-, smoothed-TV-, and Huber-TV-type ANM reconstructions, respectively. 
The red contours indicate the approximate target boundaries.

Overall, all four reconstruction methods identify the principal conductivity perturbations in most configurations. 
Their reconstruction characteristics, however, differ in terms of spatial smoothness, background homogeneity, geometric representation, and the ability to distinguish nearby inclusions.

\subsubsection{Single inclusions}
\label{subsubsec:single-inclusions}

\begin{figure}[H]
    \centering
    \includegraphics[width=0.95\textwidth, height=0.88\textheight, keepaspectratio]{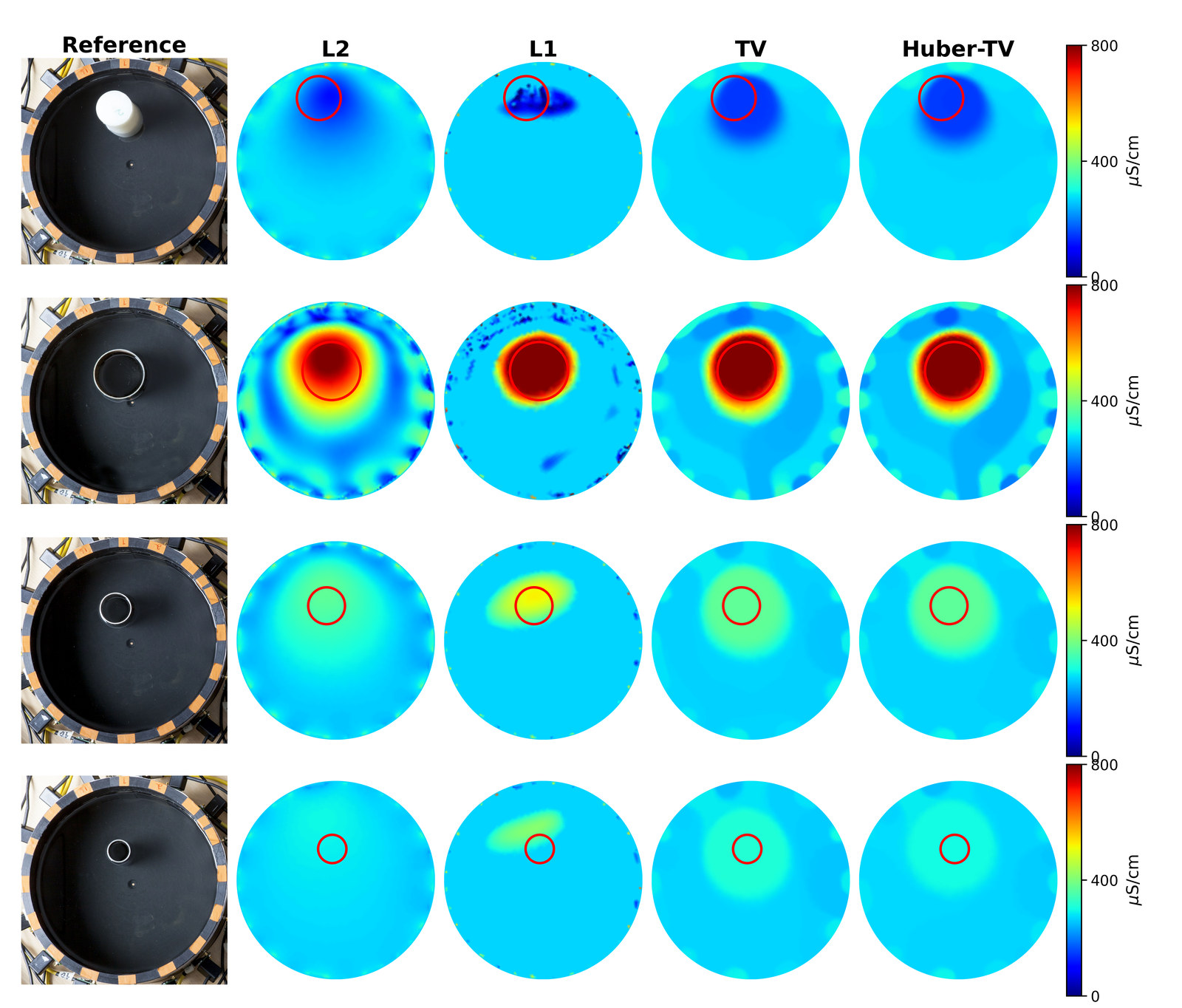}
    \caption{Single-inclusion reconstructions. }
    \label{fig:single-inclusions}
\end{figure}

Fig.~\ref{fig:single-inclusions} compares representative single-inclusion reconstructions obtained with the \(L^2\)-type, \(L^1\)-type, smoothed TV, and Huber-TV penalties for both conductive and insulating targets.

The \(L^2\)-type penalty correctly identifies the approximate inclusion locations in all four cases, but the recovered conductivity variations are distributed over relatively broad regions. 
As a result, the reconstructed interfaces are diffuse and the target extents are less accurately resolved, especially for the smaller inclusions.

The \(L^1\)-type penalty produces the most homogeneous background by penalizing deviations from the calibrated background conductivity \(\sigma_{\rm bk}\) and strongly suppressing weak conductivity variations.
This background-centered sparsity prior is particularly effective in reducing spurious variations away from the targets. 
However, it may also introduce noticeable geometric distortion. 
The reconstructed supports may become elongated, compressed, or otherwise deformed, rather than preserving the approximately circular target geometry. 
This behavior is particularly evident for the boundary-adjacent insulating target in the first row, where a localized irregular response appears near the inclusion, and for the smallest centrally located target in the final row, which is only weakly recovered.

Both TV-type penalties provide more accurate localization and sharper target--background transitions than the \(L^2\)- and \(L^1\)-type penalties.
The reconstructed target extents are generally closer to the reference configurations, although the contrast decreases and the recovered region becomes more diffuse as the inclusion becomes smaller and more distant from the boundary electrodes. 
This loss of sensitivity is most apparent for the smallest centrally located inclusion.
The smoothed TV and Huber-TV reconstructions are visually similar in these single-inclusion cases. 
Both preserve the inclusion locations and relative sizes effectively, while Huber-TV produces slightly smoother spatial transitions in some cases. 
Overall, the TV-type penalties provide the best balance among target localization, geometric fidelity, and background regularity for the single-inclusion configurations.

\subsubsection{Two inclusions}
\label{subsubsec:two-inclusions}

\begin{figure}[H]
    \centering
    \includegraphics[
        width=0.95\textwidth,
        height=0.88\textheight,
        keepaspectratio
    ]{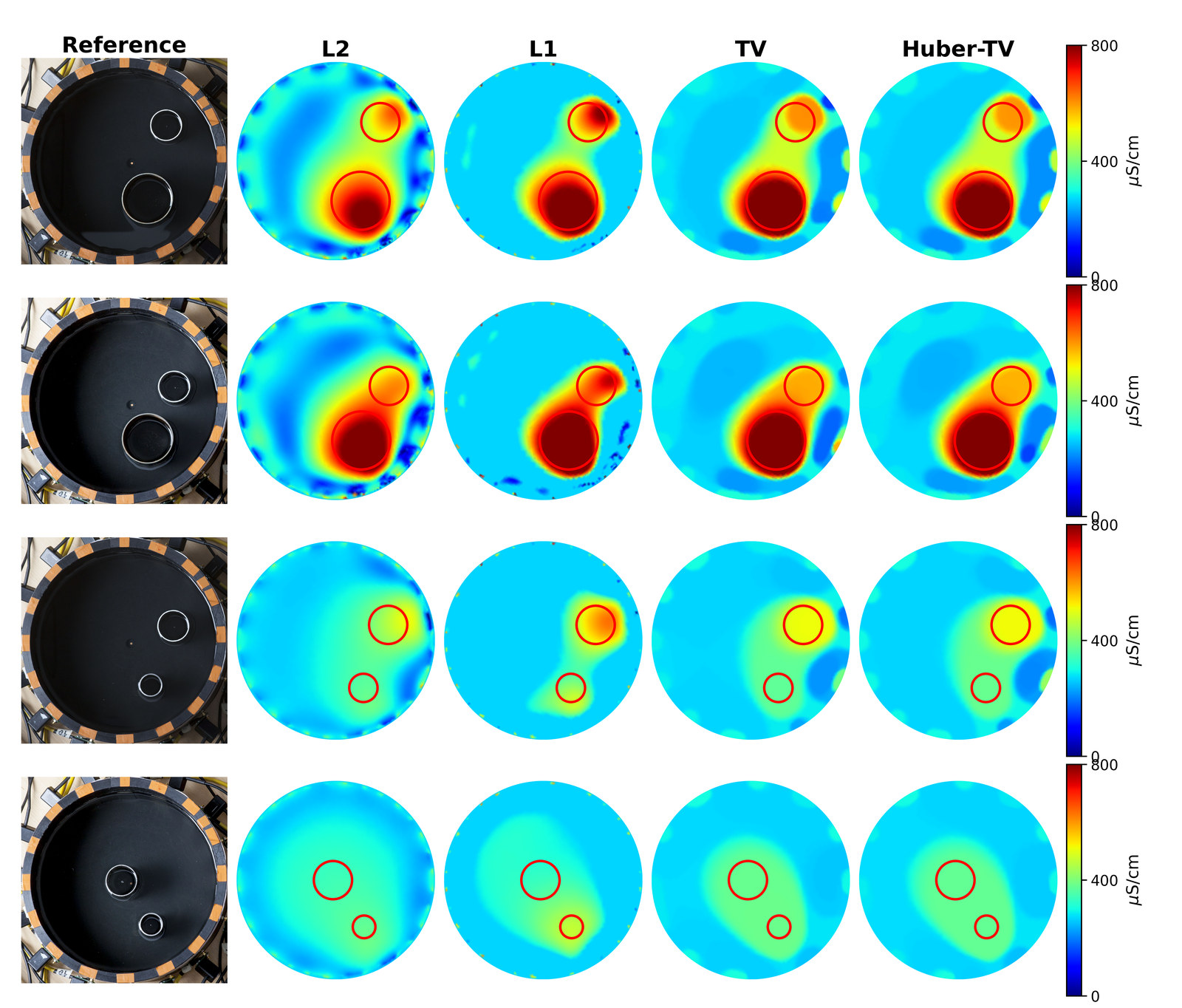}
    \caption{Two-inclusion reconstructions. }
    \label{fig:two-inclusions}
\end{figure}

Fig.~\ref{fig:two-inclusions} compares representative two-inclusion reconstructions for conductive targets with different sizes, separations, and distances from the boundary electrodes.

For the first three configurations, the \(L^2\)-type penalty identifies the approximate locations of both inclusions, but the reconstructed conductivity is distributed over broad regions. 
This produces diffuse interfaces and substantial contrast spreading, especially around the larger targets.

The \(L^1\)-type penalty suppresses weak background variations more strongly and therefore yields cleaner backgrounds. 
However, the reconstructed supports are geometrically distorted and tend to elongate toward one another. 
In the first two rows, this produces an artificial bridge between the two conductive responses and reduces their apparent separation.

The smoothed TV and Huber-TV penalties provide the most accurate localization in the first three configurations. 
They better preserve the relative target sizes and keep the two inclusions more clearly separated than the \(L^2\)- and \(L^1\)-type penalties. 
The two TV-type reconstructions are visually similar.

The final configuration is substantially more challenging because both targets are small and located farther from the boundary electrodes. 
In this case, none of the four penalties resolves the two inclusions satisfactorily.
The reconstructed responses are weak and diffuse, and the individual target sizes and boundaries cannot be reliably identified. 
This result illustrates the limited interior sensitivity and spatial resolution of EIT for small, closely spaced targets away from the electrodes.

Overall, the TV-type penalties provide the most coherent reconstructions in the less challenging cases, whereas the \(L^2\)-type penalty is affected by contrast diffusion and the \(L^1\)-type penalty by geometric distortion and artificial bridging. Nevertheless, all four methods lose resolving power for the small interior targets in the final configuration.

\subsubsection{Multiple and mixed inclusions}
\label{subsubsec:multiple-inclusions}

\begin{figure}[H]
    \centering
    \includegraphics[
        width=0.95\textwidth,
        height=0.88\textheight,
        keepaspectratio
    ]{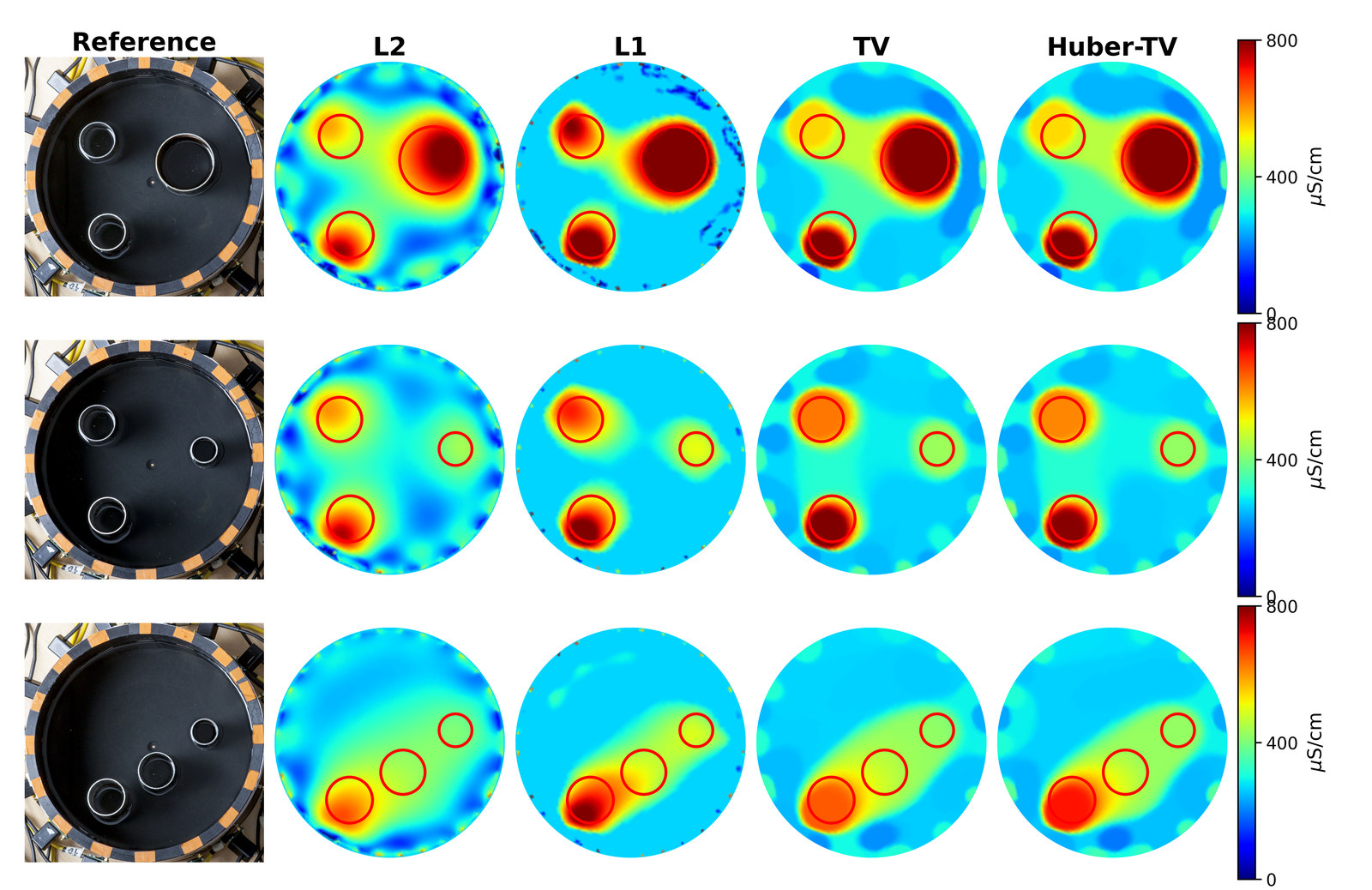}
    \caption{Reconstructions of three conductive inclusions. }
    \label{fig:three-conductive-inclusions}
\end{figure}

\begin{figure}[H]
    \centering
    \includegraphics[
        width=0.95\textwidth,
        height=0.92\textheight,
        keepaspectratio
    ]{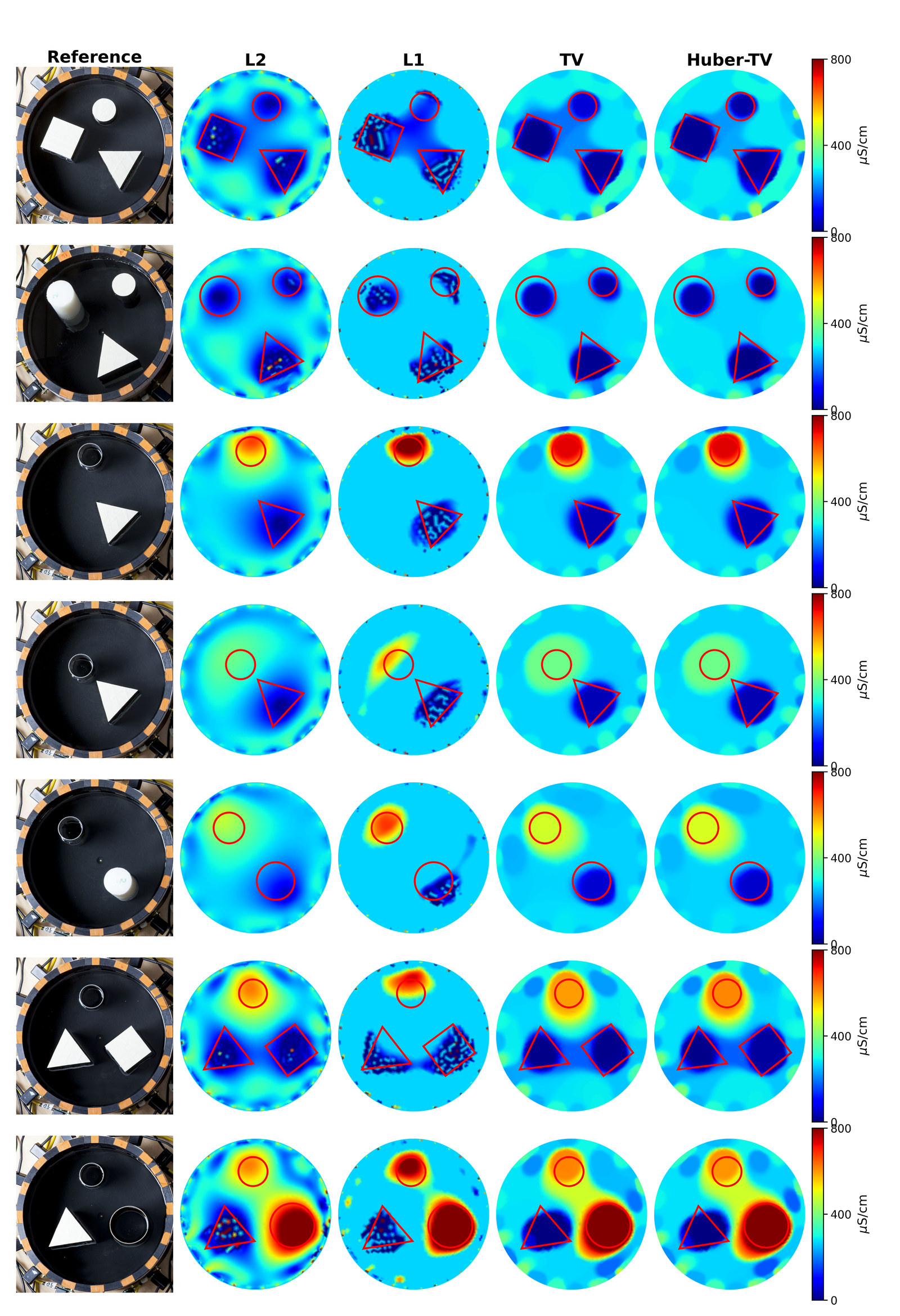}
    \caption{Reconstructions of three insulating inclusions and mixed
    conductive--insulating configurations.}
    \label{fig:mixed-inclusions}
\end{figure}

Figs.~\ref{fig:three-conductive-inclusions} and \ref{fig:mixed-inclusions} present more challenging configurations with multiple inclusions. 
Fig.~\ref{fig:three-conductive-inclusions} considers three conductive inclusions of different sizes and spatial arrangements, whereas Fig.~\ref{fig:mixed-inclusions} includes both all-insulating configurations and mixed conductive--insulating configurations.

For the three conductive inclusions in Fig.~\ref{fig:three-conductive-inclusions}, all four penalties recover the principal high-conductivity anomalies in the first two rows. 
The \(L^2\)-type reconstructions identify the approximate target locations, but the conductivity responses are spatially diffuse and accompanied by substantial background variation. 
The \(L^1\)-type penalty produces a cleaner background and separates the three targets more clearly, although the reconstructed supports are enlarged or geometrically distorted.

The smoothed TV and Huber-TV penalties provide more coherent localization and better preserve the relative arrangement of the three targets. 
Their reconstructed supports agree more closely with the reference boundaries, although some background plateaus and boundary artifacts remain. 
The third row is more difficult because the three inclusions are smaller and arranged along a narrow diagonal region. 
In this case, the individual responses overlap strongly for all four methods. 
The targets are detected only as an elongated high-conductivity structure, and their individual sizes and
boundaries cannot be recovered reliably. 
This behavior is consistent with the loss of resolving power observed for nearby targets with conductivity contrasts of the same sign.

Fig.~\ref{fig:mixed-inclusions} contains all-insulating configurations in the first two rows, followed by mixed conductive--insulating configurations.
For the all-insulating cases, all four methods identify the principal low-conductivity regions. The \(L^2\)-type reconstructions are relatively diffuse, while the \(L^1\)-type penalty produces a more homogeneous background but may yield fragmented or geometrically distorted supports. 
The two TV-type penalties provide more coherent localization, particularly for the polygonal inclusions.

In the mixed configurations, the conductive and insulating targets remain distinguishable because their conductivity contrasts have opposite signs.
The \(L^2\)-type penalty identifies the principal contrast regions but produces broader responses, whereas the \(L^1\)-type penalty yields cleaner backgrounds with occasional fragmentation of the insulating targets. 
The smoothed TV and Huber-TV reconstructions provide the most consistent separation of the opposite-sign anomalies and are visually similar across these configurations. 
Their reconstructed polygonal boundaries are smoother than the corresponding reference shapes, and some boundary artifacts remain in the more complicated cases.

\subsubsection{Challenging phantom configurations}
\label{subsubsec:challenging-phantoms}

\begin{figure}[H]
    \centering
    \includegraphics[
        width=0.95\textwidth,
        height=0.88\textheight,
        keepaspectratio
    ]{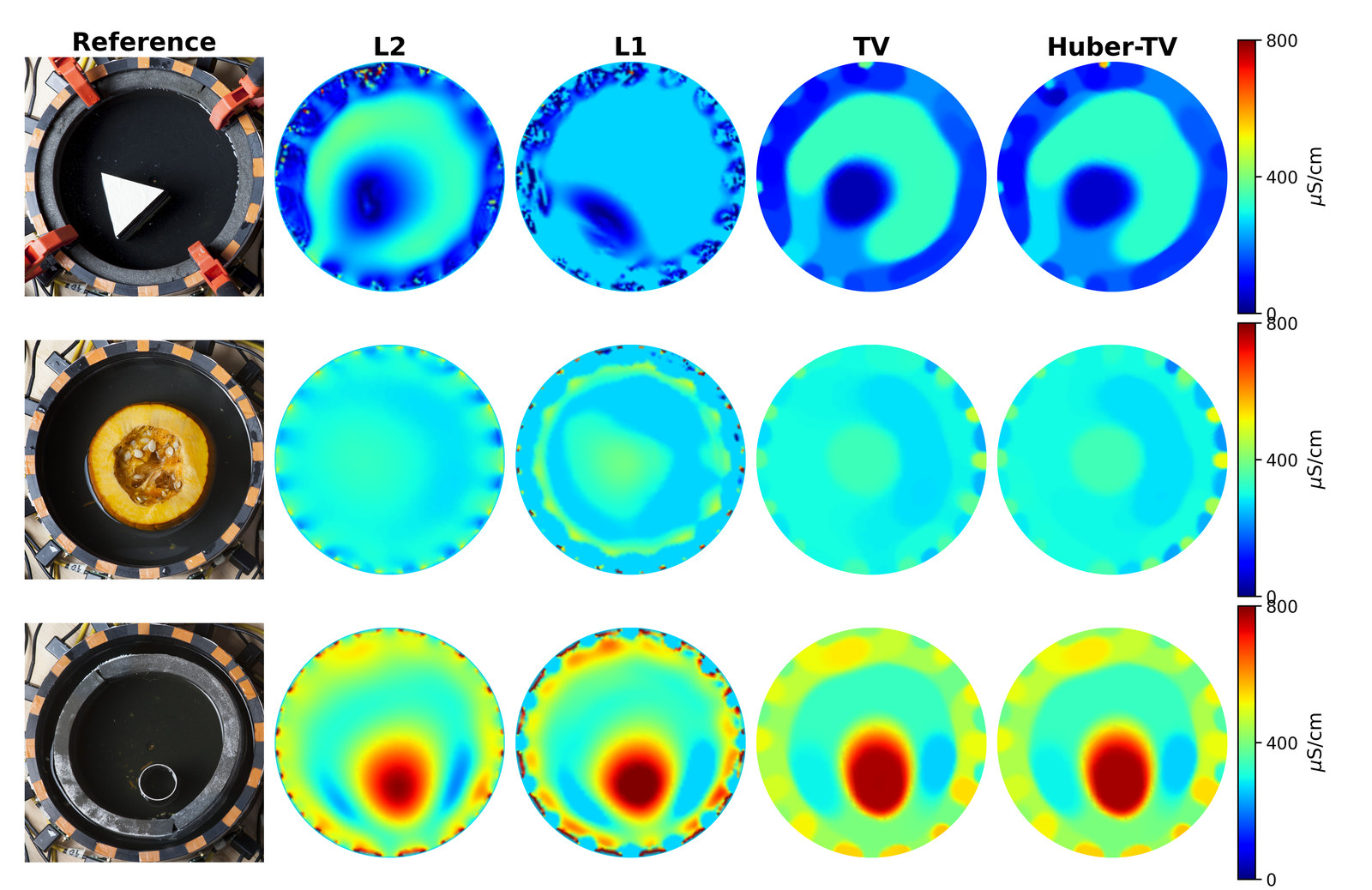}
    \caption{Reconstructions of challenging phantom configurations. }
    \label{fig:challenging-phantoms}
\end{figure}

Fig.~\ref{fig:challenging-phantoms} presents three additional challenging
phantoms: a triangular insulating target, a pumpkin, and a hollow metallic
cylinder.
For the triangular target, all methods identify a localized low-conductivity
region near the correct position. The \(L^2\)-type result is more diffuse,
whereas the \(L^1\)-type and TV-type penalties provide sharper localization,
although the reconstructed boundary is smoother than the physical triangle.
For the pumpkin phantom, all methods recover a weak conductivity perturbation
near the correct location. The \(L^1\)-type result is slightly more localized,
while the smoothed TV and Huber-TV reconstructions are smoother and visually
similar.
For the hollow metallic cylinder, all methods detect a pronounced conductive
anomaly at the correct location. The \(L^2\)-type reconstruction is more
diffuse, whereas the \(L^1\)-type and TV-type penalties produce more compact responses.
Overall, the \(L^1\)-type and two TV-type penalties show comparable
localization in these challenging examples, while the \(L^2\)-type penalty
produces smoother but more spatially diffuse reconstructions.

\subsection{Summary of numerical findings}

The numerical experiments demonstrate that the adaptive Nesterov momentum framework can be applied to measured CEM-EIT data with different structural penalties. 
The \(L^2\)-type penalty provides a smooth baseline but generally produces more spatially diffuse conductivity distributions. 
The background-centered \(L^1\)-type penalty effectively suppresses weak variations away from the targets, although the reconstructed supports may exhibit geometric deformation or fragmentation. 
The smoothed TV and Huber-TV penalties provide more coherent localization and better preserve the relative arrangement of multiple inclusions. 
Their reconstruction behavior is broadly similar across the tested KIT4 configurations.
The experiments also highlight the effect of target size, position, spacing, and conductivity contrast on the reconstruction. 
Large or boundary-adjacent inclusions are generally identified more clearly, whereas small interior targets and closely spaced inclusions of the same contrast sign produce weaker or overlapping responses.
Conductive and insulating targets in the mixed configurations remain more readily distinguishable because their conductivity perturbations have opposite signs.

\section*{Acknowledgements}

The authors thank Dr. Qiang Zhang and Dr. Dong Wang for helpful discussions on the regularity of the CEM forward and adjoint problems.

\bibliographystyle{siamplain_modified}
\bibliography{sisc_bib}

\clearpage
\appendix

\renewcommand{\thefigure}{\thesection.\arabic{figure}}
\setcounter{figure}{0}

\section{Parameter Selection for the Strong Convexity Weight}
\label{app:beta-selection}

The parameter \(\beta\) controls the relative contribution of the quadratic strong-convexity term and the structural penalty in the dual-to-primal minimization. 
Its choice therefore affects background suppression, conductivity contrast, and the geometry of the reconstructed supports.

\begin{figure}[p]
\centering
\includegraphics[width=\textwidth]
{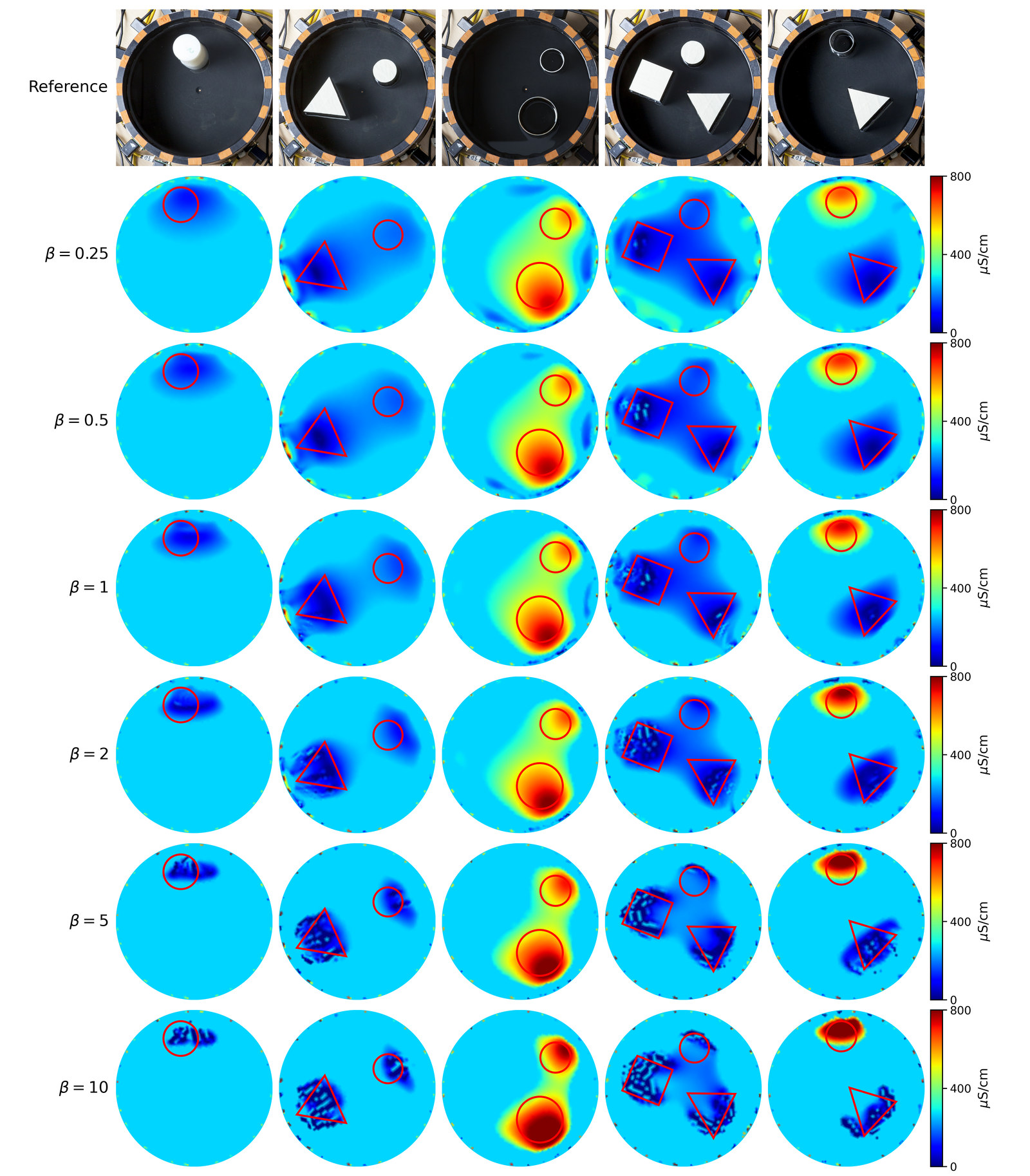}
\caption{Effect of the strong-convexity weight \(\beta\) on the
\(L^1\)-type reconstructions.}
\label{fig:beta-l1}
\end{figure}

\begin{figure}[p]
\centering
\includegraphics[width=\textwidth]
{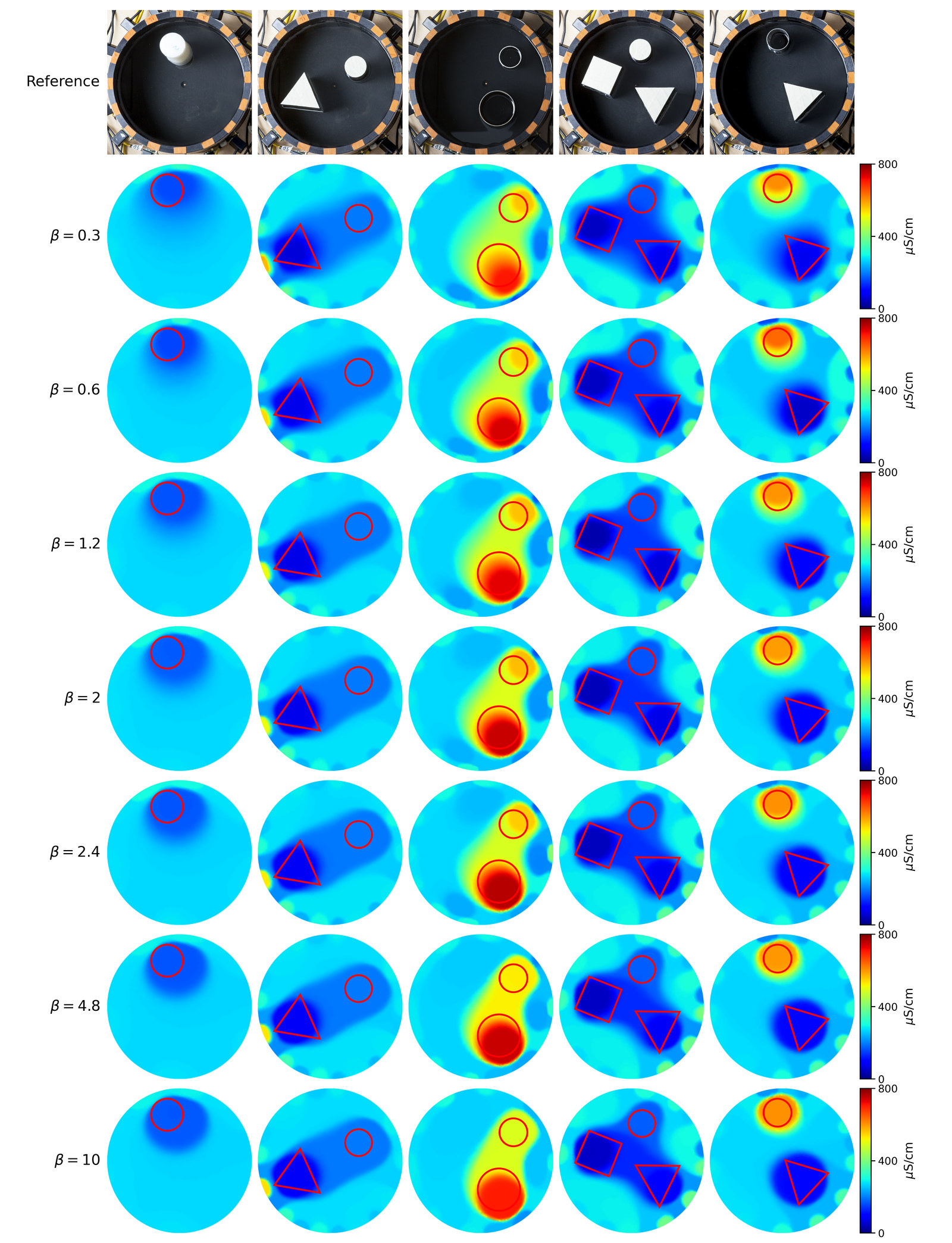}
\caption{Effect of the strong-convexity weight \(\beta\) on the smoothed TV
reconstructions.}
\label{fig:beta-tv}
\end{figure}

\begin{figure}[p]
\centering
\includegraphics[width=\textwidth]
{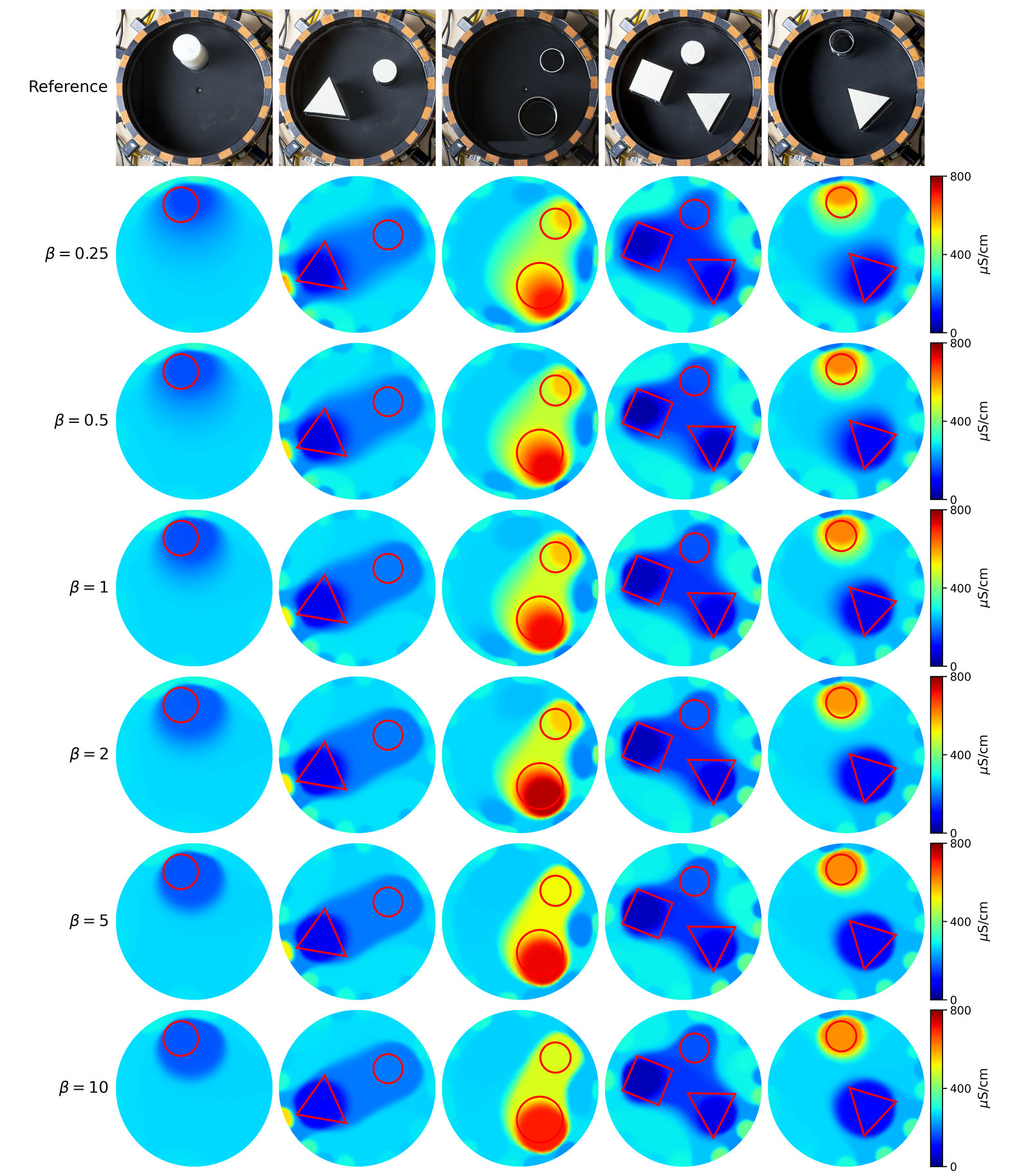}
\caption{Effect of the strong-convexity weight \(\beta\) on the Huber-TV
reconstructions.}
\label{fig:beta-huber}
\end{figure}

Since no theoretical parameter-choice rule is currently available for the measured CEM-EIT setting considered here, \(\beta\) is selected empirically through a preliminary parameter study. 
The same representative KIT4 configurations are reconstructed over a range of candidate values. 
For the
\(L^1\)-type penalty, we test
$\beta\in\{0.25,\,0.5,\,1,\,2,\,5,\,10\}$,
for the smoothed TV penalty,
$\beta\in\{0.3,\,0.6,\,1.2,\,2,\,2.4,\,4.8,\,10\}$,
and for the Huber-TV penalty,
$\beta\in\{0.25,\,0.5,\,1,\,2,\,5,\,10\}$.
The corresponding reconstructions are shown in Figs.~\ref{fig:beta-l1}--\ref{fig:beta-huber}.

For the \(L^1\)-type penalty, the effect of \(\beta\) is most evident in the background and in the spatial extent of the reconstructed anomalies. 
At the smaller values, the conductivity perturbations are relatively broad and weak background variations remain visible. 
Increasing \(\beta\) strengthens the background-centered sparsity effect and produces more localized responses.
However, the reconstructed supports also become increasingly concentrated and may exhibit fragmentation or geometric deformation, particularly for the polygonal insulating targets. 
The largest tested value, \(\beta=10\), provides strong background suppression but produces more pronounced support fragmentation in several configurations. 
The choice \(\beta_{L^1}=5\) therefore represents a compromise between background homogeneity, target contrast, and geometric representation.

The smoothed TV reconstructions are comparatively less sensitive to \(\beta\) over the tested range. 
All candidate values recover the principal conductivity perturbations and preserve their relative locations. 
Smaller values produce somewhat broader transitions and greater background variation, while intermediate values yield more coherent target regions.
At the largest values, the reconstructions remain stable, but further increasing \(\beta\) provides little visible improvement and may reduce some geometric detail. 
Among the tested values, \(\beta_{\mathrm{TV}_\varepsilon}=2\) provides a suitable balance between target localization, contrast, and background regularity across the representative configurations.

A similar degree of robustness is observed for the Huber-TV penalty. 
The principal conductive and insulating anomalies remain identifiable throughout the tested range, and the overall reconstruction pattern changes only gradually with \(\beta\). 
Increasing \(\beta\) generally sharpens the dominant responses and suppresses weak background variations, whereas excessively large values offer limited additional benefit and may emphasize local contrast at the expense of geometric smoothness. 
The value \(\beta_{\mathrm{TV}_\gamma}=2\) is selected as a moderate choice that gives consistent localization and conductivity contrast over all representative phantoms.

Overall, the parameter study shows that the \(L^1\)-type reconstruction is more sensitive to the choice of \(\beta\), especially with respect to support fragmentation and shrinkage. 
By contrast, the smoothed TV and Huber-TV formulations exhibit relatively stable reconstruction behavior over a broader range of parameter values. 
Based on the qualitative balance among background suppression, target contrast, and geometric fidelity, we use
\[
\beta_{L^1}=5,\qquad
\beta_{\mathrm{TV}_\varepsilon}=2,\qquad
\beta_{\mathrm{TV}_\gamma}=2
\]
in all numerical experiments reported in the main text.

\end{document}